\documentclass[12pt]{article}
\usepackage[dvips]{graphicx}
\begin{document}
\newcommand{\intl}{\int\limits}
\newcommand{\suml}{\sum\limits}
\setcounter{page}{1}
\tolerance=10000
\def\theequation{\arabic{section}.\arabic{equation}}
\font\title=cmbx12 scaled\magstep2 \font\bfs=cmbx12 scaled\magstep1
\font\sc=cmcsc10
\def\pni{\par\noindent}
\def\vsh{\smallskip}
\def\vs{\medskip}
\def\vvs{\bigskip}
\def\vvvs{\bigskip\medskip} 
\def\vsp{\vsh\pni} 
\def\vsn{\vsh\pni}
\def\cen{\centerline}
\def\ra{\item{a)\ }} \def\rb{\item{b)\ }}   \def\rc{\item{c)\ }}
\def\q{\quad} \def\qq{\qquad}
\def \rec#1{{1\over{#1}}}
\def\ds{\displaystyle}
\def\eg{{\it e.g.}\ }
\def\ie{{\it i.e.}\ }
\def\versus{{\it vs.}\ }
\def\e{{\rm e}}
\def\d{\partial}
\def\dx{\partial x}    \def\dt{\partial t}
\def\Ai{{\rm Ai}\,}
\def\Erfc{{\rm Erfc}\,}
\def\u{\widetilde{u}}
\def\ul{\widetilde{u}} 
\def\uf{\widehat{u}} 
\def\r{\right} \def\l{\left}
\def\rt{\right} \def\lt{\left}
\def\lra{\Longleftrightarrow}
\def\RR{\vbox {\hbox to 8.9pt {I\hskip-2.1pt R\hfil}}}
\def\NN{{\rm I\hskip-2pt N}}
\def\CC{{\rm C\hskip-4.8pt \vrule height 6pt width 12000sp\hskip 5pt}}
\def\L{{\mathcal L}} 
\def\Gz{\Gamma(z)}    \def\Ga{\Gamma(\alpha)}
\def\Gaz{\Gamma(\alpha\,,\, z)}
\def\gaz{\gamma(\alpha\,,\, z)}
\def\DG{{{D}}_\Gamma}
 \def\E{{\rm E}}
\def\Ai{{\rm Ai}\,}
\def\Erfc{{\rm Erfc}\,}
\def\Ei{{\rm Ei}\,}
\def\Ein{{\rm Ein}\,}
\def\log{{\rm log}\,}
\def\EE{{\mathcal E}}
\def\Gz{\Gamma(z)}    \def\Ga{\Gamma(\alpha)}
\def\Gaz{\Gamma(\alpha\,,\, z)}
\def\gaz{\gamma(\alpha\,,\, z)}
\def\DG{{{D}}_\Gamma}
 \def\E{{\rm E}}
\def\Ai{{\rm Ai}\,}
\def\Erfc{{\rm Erfc}\,}
\def\Ei{{\rm Ei}\,}
\def\Ein{{\rm Ein}\,}
\def\log{{\rm log}\,}
\def\EE{{\mathcal E}}
 \def\bar{\widetilde}
\def\hatt{\widehat}
\def\epsilons{{\widetilde \epsilon(s)}}
\def\sigmas{{\widetilde \sigma (s)}}
\def\fs{{\widetilde f(s)}}
\def\Js{{\widetilde J(s)}}
\def\Gs{{\widetilde G(s)}}
\def\Fs{{\widetilde F(s)}}
 \def\Ls{{\widetilde L(s)}}
 \def\Hs{{\widetilde H(s)}}
 \def\Ks{{\widetilde K(s)}}
\def\L{{\mathcal L}} 
\def\F{{\mathcal F}} 
\def\M{{\mathcal M}} 
\def\P{{\mathcal{P}}} 
\def\H{{\mathcal{H}}} 
\def\NN{{\bf N}}
\def\RR{{\bf R}}
\def\CC{{\bf C}}
\def\ZZ{{\bf Z}} 
\def\I{{\cal I}}  
\def\D{{\cal D}}  
\def\erf{\hbox{erf}}  \def\erfc{\hbox{erfc}}
\def\dfrac#1#2{\displaystyle{\frac {#1}{#2}}}

\cen{{\bf FRACALMO PRE-PRINT: \    www.fracalmo.org}}
\vsh
\cen{{\bf 
  Mathematical Methods in Economics and Finance: Vol. 1 (2006) 37-55}}
\vsh
\cen{{\sc Universit\`a Ca' Foscari di Venezia - Dipartimento di Matematica Applicata}} 
\vsh
\cen{{\sc  ISSN: 1971-6419 (Print Edition). ISSN: 1971-3878 (Electronic Edition)}}
 \vsh
 \cen{{\bf URL: http://www.dma.unive.it/mmef/}}
\vs
\hrule
\vskip 1.0truecm
\font\title=cmbx12 scaled\magstep2
\font\bfs=cmbx12 scaled\magstep1
\font\little=cmr10


\begin{center}

\bfs{The origin of infinitely divisible distributions:\\  [0.20truecm]
 from de Finetti's problem    to L\'evy-Khintchine formula}


\vvs
{Francesco MAINARDI} $^{(1)}$ and
{Sergei ROGOSIN}$^{(2)}$


\vs

$\null^{(1)}$
 {\little Department of Physics, University of Bologna, and INFN,} \\
{\little Via Irnerio 46, I-40126 Bologna, Italy} \\
{\little Corresponding Author. E-mail: {\tt francesco.mainardi@unibo.it}} 
\\ [0.25 truecm]
$\null^{(2)}$
 {\little Department of Mathematics and Mechanics, Belarusian State University,} \\
{\little Nezavisimosti Ave 4, BY-220030, Minsk, Belarus} \\
{\little E-mail: {\tt rogosin@bsu.by}}
\\ [0.25 truecm]

\end{center}

 \def\date#1{\gdef\@date{#1}} \def\@date{\today}

\vvs

\begin{abstract}
The article provides an historical survey of the
early contributions on infinitely divisible distributions
starting from the pioneering  works of de Finetti in 1929
up to the canonical forms developed in the thirties
by Kolmogorov, L\'evy and Khintchine.
Particular attention is paid to single out
the  personal contributions of the above authors
that were published in 
Italian, French or Russian
during the period 1929-1938.
In Appendix we report the translation from the Russian into English 
of a fundamental paper by Khintchine published in Moscow in 1937. 

\vs

\noindent {\bf Keywords.}
Characteristic function, infinitely divisible distributions, 
sto\-chastic processes with independent increments,
de Finetti, Kolmogorov, L\'evy, Khintchine, Gnedenko.

\vs

\noindent {\bf M.S.C. classification.}
\textsc{\small
60E07, 
60E10, 
60G51,  
01A70.}  

\vs 
\noindent {\bf J.E.L. classification.} \textsc{\small
C10, 
C16.} 
\end{abstract}


\def\pni{\par \noindent}
\def\vsh{\smallskip}
\def\s{\smallskip}
\def\vs{\medskip}
\def\vvs{\bigskip}
\def\vvvs{\bigskip\medskip} 
\def\vsp{\par}
\def\vsn{\vsh\pni}
\def\cen{\centerline}
\def\ra{\item{a)\ }} \def\rb{\item{b)\ }}   \def\rc{\item{c)\ }}
\def\eg{{\it e.g.}\ } \def\ie{{\it i.e.}\ }


\def\sg{\hbox{sign}\,}
\def\sgn{\hbox{sign}\,}
\def\sign{\hbox{sign}\,}
\def\e{{\rm e}}
\def\exp{{\rm exp}}
\def\ds{\displaystyle}
\def\dis{\displaystyle}
\def\q{\quad}    \def\qq{\qquad}
\def\lan{\langle}\def\ran{\rangle}
\def\l{\left} \def\r{\right}
\def\lt{\left} \def\rt{\right}
\def\lra{\Longleftrightarrow}
\def\d{\partial}
\def\dr{\partial r}  \def\dt{\partial t}
\def\dx{\partial x}   \def\dy{\partial y}  \def\dz{\partial z}
\def\rec#1{{1\over{#1}}}
\def\zr{z^{-1}}



\def\hatt{\widehat}
\def\epsilons{{\widetilde \epsilon(s)}}
\def\sigmas{{\widetilde \sigma (s)}}
\def\fs{{\widetilde f(s)}}
\def\Js{{\widetilde J(s)}}
\def\Gs{{\widetilde G(s)}}
\def\Fs{{\wiidetilde F(s)}}
 \def\Ls{{\widetilde L(s)}}
\def\L{{\mathcal L}} 
\def\F{{\mathcal F}} 


\def\NN{{\rm I\hskip-2pt N}}
\def\RR{{\bf R}}
\def\CC{{\bf C}}
\def\ZZ{{\bf Z}} 


\def\I{{\cal I}}  
\def\D{{\cal D}}  

\def\erf{\hbox{erf}}  \def\erfc{\hbox{erfc}}
\def\intR{\int_{-\infty}^{+\infty}}

 \def\EE{\vbox {\hbox to 8.9pt {I\hskip-2.1pt E\hfil}}}
 \def\PP{\vbox {\hbox to 8.9pt {I\hskip-2.1pt P\hfil}}}
\def\GG{{\mathcal{G}}}


\newpage
 
\section{Introduction}

The purpose  of this paper is
to illustrate  how the concept of an infinitely divisible distribution
has been developed up to obtain the canonical form of its
characteristic function.

Usually  historical aspects  on this development are
known thanks to  some notes  available in the classical textbooks
by L\'evy \cite{Levy ADDITION37}
(published in French in 1937 and 1954),
  by Gnedenko-Kolmogorov \cite{Gnedenko-Kolmogorov LIMITDISTR54}
(published in Russian in 1949
and translated into English in 1954)
and by Feller \cite{Feller 2}
(published in English in 1966 and  1971).
Similar historical notes can be extracted 
from the recent treatises by Sato \cite {Sato BOOK99} 
and by Steutel and van Harn \cite{Steutel-vanHarn BOOK04}.

In our opinion, however, a better  historical  analysis can be accomplished if
one examines the original works of the pioneers, namely
Bruno de Finetti (1906-1985)
\cite{DeFinetti 29a,DeFinetti 29b,DeFinetti 29c,DeFinetti 30h,DeFinetti 31c},
 Andrei Nikolaevich Kolmogorov (1903-1987)
\cite{Kolmogorov 32a,Kolmogorov 32b},
who published in Italian in  { \it Rendiconti della R. Accademia Nazionale dei Lincei},
 Paul L\'evy (1886-1971) \cite{Levy 34,Levy 35},
who published in French  in   {\it Annali della  R. Scuola Normale di  Pisa},
and finally  Alexander Yakovlevich Khintchine
\footnote{There is also the transliteration  Khinchin.}
(1894-1959) \cite{Khintchine 37},
who published in Russian in the
 { \it Bulletin of the Moscow State University}.
 Noteworthy is the  1938 book by Khintchine himself
\cite{Khintchine LIMITDISTR38}, in Russian,
on {\it Limit Distributions for Sums
of Independent Random Variables}.
For the  reader interested in  the biographical notes and bibliography
of the mentioned scientists 
we refer: 
for de Finetti to \cite{Cifarelli-Regazzini 96,Daboni UMI87,Daboni ATTI87,DeFinetti OPERE06},
for Kolmogorov to \cite{Shiryaev 92,Shiryaev 00},
for L\'evy  to \cite{Levy LIFE70,Levy COLLECTEDPAPERS73}
and for Khintchine to \cite{Gnedenko OBITUARY}.

In spite of the fact that de Finetti
was the pioneer of  the infinitely divisible distributions
in view of his 1929-1931 papers, as is well recognized in the literature,
the attribute {\it infinitely divisible}, as noted by  Khintchine
 in his 1938 book
\cite{Khintchine LIMITDISTR38}, first appeared in the
Moscow mathematical school, precisely
in the 1936 unpublished thesis by G.M. Bawly (1908-1941)
\footnote{Gregory Minkelevich Bawly 
graduated at the Moscow State
University in 1930, defended his PhD thesis under guidance of
A.N. Kolmogorov in 1936. His scientific advisor had greatly
esteemed his results on the limit distributions for sums of
independent random variables and cited him in his book
with Gnedenko \cite{Gnedenko-Kolmogorov LIMITDISTR54}.
G.M. Bawly lost his life  in
Moscow in November 1941 at a bombing attack.
}.
According to Khintchine  \cite{Khintchine LIMITDISTR38}
the  name of  {\it infinitely divisible
distributions} (in  a printed version)  is  found
in the 1936 article by G.M. Bawly \cite{Bawly 36},
that was recommended
for publication in the very important starting volume of the new
series of Matemati\v{c}eski Sbornik.
However, we note 
that this term was not ``stably'' applied in the article. Two alternative
(and equivalent) terms were used, namely  infinitely =
unbeschr\"{a}nkt (German) and unboundedly = unbegrenzt (German),
see \cite[p. 918]{Bawly 36}.

The first formal
definition of an infinitely divisible distribution was given by
Khintchine himself \cite{Khintchine 37a}. It reads: 
{\it a distribution of a random variable which for any positive integer
$n$ can be represented as a sum of $n$ identically distributed
independent random variables is called an infinitely divisible
distribution}.

We note that infinitely
divisible distributions (already under this name) were formerly
studied systematically in the 1937 book by L\'{e}vy \cite{Levy ADDITION37},
and soon later in the 1938 book by Khintchine \cite{Khintchine LIMITDISTR38}.
We also note that  L\'evy himself, in his  late biographical
1970 booklet \cite[p. 103]{Levy LIFE70},  attributes to Khintchine
the name {\it ind\'efiniment divisible}.
The canonical form of infinitely divisible
distributions is known in the literature as {\it L\'evy-Khintchine formula},
surely because it was so named by  Gnedenko and Kolmogorov
\cite{Gnedenko-Kolmogorov LIMITDISTR54}
in their
classical treatise on {\it Limit Distributions for Sums of Independent
Random Variables}
\footnote{We note that the Russian titles of both books by
Khintchine and Gnedenko \& Kolmogorov are identical, although in
the reference list (in Russian and in English) of the book by
Gnedenko \& Kolmogorov the title of Khintchine's previous book is
in some way different ({\it Limit Theorems  for Sums of Independent
Random Variables}).} that has appeared in Russian in 1949 and in English in 1954.

The plan of  the present    paper is as follows.
In Section 2 we provide a survey of the known results
on infinitely divisible distributions.
Then we pass
to  present  the tale on the
origin     of these results
by recalling, in a historical perspective, the early publications
of our four actors:
 de Finetti, Kolmogorov, L\'evy and Khintchine.
Section 3 is  devoted to de Finetti and Kolmogorov,
namely to
the so-called de Finetti's problem (as it was referred to by Kolmogorov).
Section 4 is devoted to L\'evy and   Khintchine,
namely to the origin of the so-called L\'evy-Khintchine formula.

To our knowledge  the original contributions
by Khintchine have never been translated into English,
so we find it convenient to report in Appendix
the English translation of his 1937  paper, that
has led to the L\'evy-Khintchine formula.
We plan to publish the English translation 
of the 1938 book
by Khintchine on  {\it Limit Distributions for Sums
of Independent Random Variables}
along with
a few related articles of him (originally in Italian, German and Russian).

Concerning our bibliography, the main text and the footnotes    give
references to  some classical publications.
However, we take this occasion
to edit a more  extended
bibliography on infinite divisible distributions and related topics,
that, even if non-exhaustive, could be of some interest.

\section{A survey  on   infinitely 
 divisible \\ distributions}

Hereafter we recall  the classical results  on
infinitely divisible distributions just to introduce our notations.
We presume that the reader has a good knowledge in the Probability Theory.
In the below formulations we essentially follow the 
treatments by Feller \cite{Feller 2} and  by Lukacs \cite{Lukacs CF70};
 in the references, however,  we have cited several treatises
 containing excellent chapters on infinite divisible distributions.    

A probability distribution $F$ is infinitely divisible iff for each
$n \in \NN$ it can be represented as the distribution of the sum
$$ S_n = X_{1,n} + X_{2,n} + \ldots + X_{n,n} \eqno(2.1) $$
of $n$ independent random variables with a common distribution $F_n$.
It is common to locate  the random variables in an infinite  triangular
array
$$
\begin{array}{l}
 X_{1,1} \\
 X_{2,1},\, X_{2,2} \\
 X_{3,1},\, X_{3,2},\, X_{3,3} \\
\hspace{10mm} \ldots \\
 X_{n,1}, \, X_{n,2}, \, X_{n,3}, \,\ldots, \,X_{n,n} \\
\hspace{10mm} \ldots    \hspace{10mm} \ldots
\end{array}
\eqno(2.2)
$$
whose rows 
contain independent  identically distributed ($iid$) random variables.

This definition is valid in any number of dimensions,  but for
the present we shall limit our attention to one-dimensional
distributions. It should be noted that infinite divisibility is a property
of the {\it type}, that is, together with $F$ all distributions
differing from $F$ only by location parameters
are infinitely divisible.
{\it Stable} distributions
(henceforth the Gaussian and the Cauchy distributions)
are infinitely  divisible   and distinguished by the fact that $F_n$
differs from $F$ only by location parameters.

On account of the convolution property
of the  distribution functions of independent random variables,
the  distribution function $F$ turns out to be  the $n$-fold
convolution of some distribution function $F_n$; then,
 the notion of infinite divisibility can be introduced  by means of the
{\it characteristic function}:
$$  \varphi (t)  := \EE \{\e^{\ds \, it X}\} := \int_{ - \infty}^{+\infty}
               \e^{\ds \, it x} \, dF(x)\,.
\eqno(2.3)$$
In fact, for an infinitely divisible distribution its characteristic
function $\varphi(t)$ turns out to be, for every positive integer $n$,
 the $n$-th power of some characteristic function. This means
that there exists, for every positive integer $n$, a characteristic function
$\varphi _n(t)$ such that
$$  \varphi  (t)  = [\varphi _n(t)]^n\,. \eqno(2.4)$$
The function  $\varphi _n(t)$ is uniquely determined by $\varphi (t)$,
$  \varphi _n(t) = [\varphi(t)]^{1/n}$, provided that
one selects  the principal branch for the $n$-th root.

Since Eqs. (2.2) and (2.4) are equivalent, alternatively   one could
speak about  infinitely divisible distributions
or  infinitely divisible  characteristic functions.
Elementary  properties  of infinitely divisible characteristic functions
are listed by Lukacs \cite{Lukacs CF70}.
The concept of infinite divisibility is very important
in probability theory, particularly in the study of limit theorems.


Here  we stress the
fact that infinitely divisible distributions are intimately connected
with {\it stochastic processes with independent increments}.
By this  we mean a   family of random variables $X(\lambda) $
depending on the continuous parameter $\lambda $ and such that
the increments $X(\lambda_{k+1}) -X(\lambda _k)$ are
mutually independent for any finite set
$\{\lambda _1<\lambda _2< \ldots  <\lambda _n\}$.
More precisely the processes are assumed  to be {\it homogeneous},
that is with {\it stationary increments}. Then
the distribution
of $Y(\lambda ):= X(\lambda _0 +\lambda )-X(\lambda _0)$
depends only on the length  $\lambda $ of the interval but not on $\lambda _0$.
Let us make a partition the interval $[\lambda _0, \lambda _0+\lambda]$
by $n+1$ equidistant points
$\lambda _0<\lambda _1<\ldots <\lambda _n =\lambda _0+\lambda $
and put $X_{k,n} =X(\lambda _k) -X(\lambda _{k-1})$.
Then the variable $Y(\lambda )$ 
of a process with stationary independent increments is the sum of
$n$ independent variables $X_{k,n}$ with a common
distribution and hence $Y(\lambda )$ 
{\it has an infinitely divisible distribution}.
We can summarise  all above by simply writing
 the characteristic function of $Y(\lambda )$  for any $\lambda >0\, $ as
$$\varphi(t,\lambda ) := \EE \left \{\e^{\ds \, it X(\lambda)}\right \} =
   \{ \varphi(t,1)\}^\lambda \,.\eqno(2.5) $$

We note that we have adopted the notation
commonly used in  the early contributions:
the letter $t$
denotes the Fourier parameter of the  characteristic function
whereas  the continuous parameter (essentially the time) of a stochastic process
has been denoted by  the letter $\lambda $.
Only later, when the theory of stochastic processes
became well developed,
the authors had denoted  the Fourier parameter
by a different letter like $u$ or $\kappa $
reserving, as natural,  the letter $t$ to the time
entering the stochastic processes.
The reader should    be aware of the old notation in order
to avoid possible confusion.

Let us close this section by recalling (essentially based on the book by Lukacs)
the  main theorems concerning the structure
of infinitely  divisible distributions,
that are relevant to our historical survey.
\vskip 0.25truecm
\noindent
{\bf First de Finetti's Theorem}:
A characteristic function is infinitely divisible iff
it has the form
$$ \varphi (t) = \lim_{m\to \infty} \exp\{ p_m[\psi_m(t) -1]\}\,, \eqno(2.6)$$
where the $p_m$ are real positive numbers while $\psi_m(t)$
are characteristic functions.
\vskip 0.25truecm
\noindent
{\bf   Second de Finetti's Theorem}:
The limit of a sequence of finite products
of Poisson-type characteristic functions is infinitely divisible.
The converse is also true.
This means that the class of infinitely divisible laws coincides with the class
of distribution limits of  finite convolutions of  distributions of Poisson-type\footnote{
Let us recall that for the characteristic function of the Poisson distribution we have
according to (2.4)
$$  \varphi(t) = \exp\left[\lambda\,\left(\e^{it}-1\right)\right] \,, \quad \hbox{so that} \quad
\varphi _n(t) = \exp\left[\frac{\lambda}{n}\,\left(\e^{it}-1\right)\right] \,.$$
The theorem can be used to show that a given characteristic function is
infinitely divisible. For an example we refer the reader to \cite[p. 113]{Lukacs CF70}.}.
\vskip 0.25truecm
\noindent
{\bf The Kolmogorov canonical representation}:
The function $\varphi (t)$ is the characteristic function
of an infinitely divisible distribution with finite second moment iff
it can be written in the  form
$$ \log \varphi (t) = i\gamma t  
+ \int_{-\infty}^{+\infty}
\left( \e^{\ds\, itu} -1 -itu \right) \, \frac{dK(u)}{u^2}\,, \eqno(2.7)$$
where $\gamma $ is a real constant,
 and  $K(x)$ is a non-decreasing and bounded function
such that $K(-\infty)=0$.
The integrand is defined for $u=0$ to be equal to $-(t^2/2)$.
\vskip 0.25truecm
\noindent
{\bf The L\'evy  canonical representation}:
The function $\varphi (t)$ is the characteristic function
of an infinitely divisible distribution  iff
it can be written in the  form
$$
\begin{array}{ll}
\log \varphi (t) = i \gamma t -{\ds \frac{\sigma ^2}{2}\,t^2}
& +
{\ds  \int_{-\infty}^{-0}
\left( \e^{\ds\, itu} -1 -\frac{itu}{1+u^2} \right) \, dM(u)} \\
& +
{\ds  \int_{+0}^{\infty}
\left( \e^{\ds\, itu} -1 -\frac{itu}{1+u^2} \right) \, dN(u)}
\,,
\end{array}
\eqno(2.8)$$
where $\gamma $ is a real constant,
$\sigma ^2$ is a real and non-negative constant,
and the functions $M(u)$, $N(u)$
satisfy the following conditions:
\\
(i) $M(u)$ and $N(u)$ are non-decreasing in   
($-\infty,0$) and ($0, +\infty$), respectively.
\\
(ii) $M(-\infty)= N(+\infty) = 0\,. $
\\
(iii) The integrals $\int_{-\epsilon }^0 u^2\, dM(u) $
   $\int_0^{+\epsilon }  u^2\, dN(u) $ are finite for every $\epsilon >0$.
\vskip 0.25truecm
\noindent
{\bf The L\'evy-Khintchine   canonical representation}:
The function $\varphi (t)$ is the characteristic function
of an infinitely divisible distribution  iff
it can be written in the  form
$$
\log\, \varphi(t) = i \gamma t
+  \int_{-\infty}^{+\infty}
\left[\e^{\ds \, i  t u} - 1 - \frac{i t u}{1 + u^2}\right]
\frac{1 +u^2}{u^2} \, d G(u)\,,
\eqno(2.9) $$
where $\gamma $ is a real constant,
 and  $G(u)$ is a non-decreasing and bounded  function
such that $G(-\infty)=0$.
The integrand is defined for $u=0$ to be equal to $-(t^2/2)$.

We point out that there is  a  tight connection between  the L\'evy-Khintchine   canonical representation 
and  the general  Central Limit Theorem. 
For a clear description of a modern view on this connection we refer, {\it e.g.},
to  \cite{Hoffman 94}.
\section{The work of de Finetti and Kolmogorov}

Bruno de Finetti is recognized to be  the most prominent scientist 
of the Italian school of Probability and Statistics,
that started at the beginning of the last century with
 Guido Castelnuovo (1865-1952) and  Francesco Paolo Cantelli (1875-1966).
His personality and his interest in probability
came out already with his  attendance
at the  1928 International Congress
of Mathematicians\footnote{
The  Chairman  of the Congress was  Salvatore Pincherle (1853-1936),
Professor of Mathematics at the
University of Bologna from 1880 up to  1928, the year of the Congress.
He was the first President  of the Unione Matematica Italiana (UMI)
from 1922 up to 1936, at his death.}
held  in  Bologna (Italy)
 from 3 to 10 September 1928.
The young de Finetti presented a note on the role
of the characteristic
function in random  phenomena \cite{DeFinetti BOLOGNA28}\footnote{We have to mention that 
this work by de Finetti was the first significant contribution to the subject known
now as the theory of exchangeable sequences (of events).
A more exhaustive account appeared in 1931 in \cite{DeFinetti MEM30}.}, 
that  was  published only
 in 1932   in the Proceedings of the Congress\footnote{
The Proceedings
were published by Zanichelli, Bologna, with all  the details of
the scientific and social programs, in 6  volumes,
that appeared from 1929 to 1932.
The papers, published in one of the following languages:
Italian, French, German and English,
were classified  in 7 sessions according to
their topic.
The papers presented by Cantelli \cite{Cantelli BOLOGNA28},
de Finetti, Romanovsky and Slutsky
(Session IV,  devoted to Actuarial Sciences, Probability and  Statistics)
were included  within the last  volume, published in 1932.}.

At the Bologna Congress de Finetti had  the  occasion  to meet
L\'evy and Khintchine  (who were included in the
French and Russian  delegations, respectively)
but we are not informed about their interaction.
We note, however, that Khintchine did not present any communication
whereas
 L\'evy  presented a note outside the field
of probability, precisely on fractional
differentiation \cite{Levy BOLOGNA28};
furthermore Kolmogorov did not attend the Congress.
Surely L\'evy, Kolmogorov and Khintchine held
in high  consideration the Italian school of Probability
since  in  the thirties
they submitted some relevant papers
to  Italian journals
(written in  Italian for the Russians and in French for L\'evy),
see \eg \cite{Khintchine ATTUARI32,Khintchine ATTUARI35,Khintchine ATTUARI36,
Kolmogorov 32a,Kolmogorov 32b,Levy 34,Levy 35}.

Just after the Bologna Congress
 de Finetti started a research regarding
functions with random increments, see
\cite{DeFinetti 29a,DeFinetti 29b,DeFinetti 29c,
DeFinetti 30h,DeFinetti 31c}
based  on the theory of infinitely divisible
characteristic functions, even if he did not use such term.
His  results can be summarized in a number of
relevant theorems   (partly stated in the previous Section). As
it was already mentioned they are highly connected with 
the stochastic processes
with stationary independent increments.
In this respect we refer the reader to the Section 2.2
of the excellent paper by
Cifarelli and Regazzini on de Finetti's contributions
in Probability and Statistics  \cite{Cifarelli-Regazzini 96}.

The papers by de Finetti,
 published  in the period 1929-1931
(in Italian) in  the Proceedings of the Royal Academy of Lincei
({\it Rendiconti della Reale Accademia Nazionale  dei Lincei})
\cite{DeFinetti 29a,DeFinetti 29b,DeFinetti 29c,
DeFinetti 30h,DeFinetti 31c},
attracted the attention of Kolmogorov
who was interested to solve  the so-called
{\bf de Finetti's problem}, that is  the problem
of finding the general formula
 for the characteristic function of
the infinitely divisible distributions.
This problem  was indeed attacked
by Kolmogorov in 1932 in two notes published in Italian
in the same journal as de Finetti (Rendiconti della Reale Accademia Nazionale  dei Lincei),
where he gave    an exhaustive answer to {\bf de Finetti's problem}
for the case of variables  with
{\it finite second moment}, see
\cite{Kolmogorov 32a,Kolmogorov 32b}.
These two notes are available in English in a unique  paper (No 13)
in the {\it Selected Works of A.N. Kolmogorov}
with a comment of  V.M. Zolotarev, see   \cite{Kolmogorov PS13}:
the final result of Kolmogorov is reported in Section 2 as Eq. (2.7),
known as {\bf the Kolmogorov canonical representation}
of the infinitely divisible characteristic functions.

\section{The work of L\'evy and  Khintchine}

 The general case of de Finetti's problem, including also the case of
{\it infinite variance}, was investigated 
in 1934-35 by L\'evy \cite{Levy 34,Levy 35}
who published  two papers in French in the Italian
Journal: Annali della Reale Scuola Normale di Pisa.
At that time L\'evy was quite interested in the so-called
{\it stable distributions} that are known to exhibit  infinite variance,
except for the particular case of the Gaussian.

The approach by L\'evy,
 well described in his classical 1937 book \cite{Levy ADDITION37},
is quite independent  from that of Kolmogorov,
as can be understood from footnotes
in his 1934 paper \cite{Levy 34},
that we report partly below in original.
From the foot\-note$^{(1)}$ we learn that the results  contained in his paper
were presented in three communications  of the Academy of Sciences
(Comptes Rendus) of 26 February, 26 March and 7 May 1934.
Then, in  the footnote    $^{(6)}$, p. 339, the Author writes:
\\
{\it [Ajout\'e \'a la correction des \'epreuves]
Le r\'esum\'e de ma note du 26 f\'evrier, r\'edig\'e
par M. Kolmogorov, a attir\'e mon attention sur deux Notes
de M. B. de Finetti} (see \cite{DeFinetti 29a,DeFinetti 30h})
{\it et deux autres de M. Kolmogorov lui-m\^eme}
(see \cite{Kolmogorov 32a,Kolmogorov 32b}),
{\it  publi\'ees dans
les Atti Accademia Naz. Lincei (VI ser).
Ces derni\`eres notamment contiennent la solution du probl\`eme
trait\'e dans le pr\'esent travail, dans le cas o\`u
le processus est homog\`ene et o\`u la valeur probable
$\EE\{x^2\}$ est finie.
Le r\'esultat fondamental du pr\'esent M\'emoire
apparait donc comme une extension d'un r\'esultat de M. Kolmogorov.}

This  means that P.~L\'{e}vy was not aware  about
the results on homogeneous processes with independent increments
obtained
by B.~de Finetti and  by A.~N.~Kolmogorov.
The final result of L\'evy is reported in Section 2 as Eq. (2.8),
known as {\bf the L\'evy canonical representation}
of the infinitely divisible characteristic functions.

In a paper of 1937 Khintchine  \cite{Khintchine 37} showed
that L\'evy' s result can be obtained also by
an extension of Kolmogorov's method:
his final  result,
reported in Section 2 as Eq. (2.9),  is known
as {\bf the L\'evy-Khintchine  canonical representation}
of the infinitely divisible characteristic functions.
The translation from the Russian of this fundamental paper
can be found in Appendix.
The theory of infinitely distributions  was then  presented  in
German in the article \cite{Khintchine MatSb37}
and in Russian in  his 1938 book on
{\it Limit Distributions  for Sums of Independent Random Variables}
\cite{Khintchine LIMITDISTR38}.

Unfortunately, many  contributions by Khintchine
(being in Russian) remained almost unknown
in the West
up to the English translation of the
treatise by Gnedenko and Kolmogorov \cite{Gnedenko-Kolmogorov LIMITDISTR54} in 1954.

The obituary of Khintchine \cite{Gnedenko OBITUARY},
that Gnedenko (his former pupil) presented at the
 1960 Berkeley Symposium
on Mathematical Statistics and Probability,
 provides a general  description of the works of Khintchine
 along with a complete bibliography.
From that we learn that the 1938 book by Khintchine
was preceded by a special course of lectures
in Moscow University that attracted the interest of
A.A. Bobrov, D.A. Raikov and B.V. Gnedenko himself.


\section*{Acknowledgements}

The authors are grateful to 
R. Gorenflo  
and the anonymous referees
for useful comments. 
We thank also
O. Celebi for the help with the paper by G.M. Bawly published in Turkey.

\section*{Appendix: Khintchine's 1937 article}

\noindent
{\bf {A.~Ya. Khintchine}}\footnote{
We have to remark that the footnotes in this
Appendix are translation of the original ones by Khintchine.}
:
A new derivation of a formula by P.~L\'evy,
\\
{\it Bulletin of the Moscow State University}
 {\bf 1} (1937) 1-5.


\vskip 0.25truecm
\noindent
A collection of all the so-called infinitely divisible distributions
was discovered for the first time by P.~L\'evy\footnote{
 Ann. R. Scuola Norm. Pisa (Ser. II), {\bf 3}, pp. 337-366 (1934).}.
He has derived a  remarkable formula for the logarithm of the
characteristic function of such a distribution. Because of the
importance of this formula I shall give here a new completely
analytic and very simple proof of it\footnote{%
The method of this proof can be considered
as an extension of the idea by A.~N.~Kolmogorov.
The latter formed the base of the proof of an
analogous formula in the important case of finite variance (see
{\it Rendiconti dei Lincei}, {\bf 15}, pp. 805-808 and 866-869
(1932).}.

Let $\varphi(x)$ be an infinitely divisible distribution and let
$\varphi(t)$ be the corresponding characteristic function. It is
known that for each $h\geq 0$ the function ${\varphi(t)}^h$ is a
characteristic function as well. We denote by $\varphi_h(x)$ the
corresponding distribution. Thus
$$
\log\, {\varphi(t)} = \lim_{h\to 0} \frac{{\varphi(t)}^h -
1}{h} = \lim_{h\to 0} I_h(t),
$$
where
$$
I_h(t) = \frac{1}{h} \, \int\limits_{-\infty}^{+\infty} \left(\e^{i t
u} - 1\right)\, d \varphi_h(u)\,.
$$
Put for each $h > 0$
$$
G_h(u) =  \int\limits_{0}^{u} \frac{v^2}{v^2 + 1}
\frac{d \varphi_h(v)}{h} \,.
\eqno(A.1) $$
Clearly the  function $G_h(u)$ is nondecreasing and bounded.
Furthermore,
$$
I_h(t) =  \int\limits_{-\infty}^{+\infty} \left(\e^{i t u} -
1\right) \frac{u^2 + 1}{u^2}\,d G_h(u)\,.
$$
Taking the real part of this formula we have
$$
-{\mathrm{Re}}\, I_h(t) =
\int\limits_{-\infty}^{+\infty} \left(1 - \cos\, t u\right)
\frac{u^2 + 1}{u^2}\, d G_h(u)\,.
\eqno(A.2)
$$
Let
$$
A_h := \int\limits_{|u|\leq 1} d G_h(u),\;\; B_h :=
\int\limits_{|u|> 1} d G_h(u),\;\; C_h :=  A_h +  B_h =
\int\limits_{-\infty}^{+\infty} \,d G_h(u)\,.
$$
Relation (A.2) gives us for $t = 1$
$$
-{\mathrm{Re}}\, I_h(1) \geq \int\limits_{|u|\leq 1} \left(1 -
\cos\, t u\right) \frac{u^2 + 1}{u^2}\, d G_h(u) \geq c A_h\,,
\eqno(A.3) $$
where $c$ is a strictly positive  constant. In the same way for
each $t$ we have
$$
-{\mathrm{Re}}\, I_h(t) \geq \int\limits_{|u|\geq 1} \left(1 -
\cos\, t u\right) \,d G_h(u)\,.
$$
Hence
$$ -\int\limits_{0}^{2} {\mathrm{Re}}\, I_h(t) d t \geq 2 B_h -
\int\limits_{|u|\geq 1} \frac{\sin\, 2 u}{u} \,d G_h(u) \geq B_h\,.
\eqno(A.4)$$
It follows from Eqs. (A.3)-(A.4) that
$$
C_h = A_h + B_h = -\frac{{\mathrm{Re}}\, I_h(1)}{c} -
\int\limits_{0}^{2} {\mathrm{Re}}\, I_h(t) \,d t\,.
$$
Since the function $I_h(t)$ uniformly converges on $0\leq t \leq 2$
as $h\to 0$ to a finite limit, then $C_h$ is bounded
as $h \to 0$. Since $G_h(0) = 0$,  the functions $G_h(u)$
remain uniformly bounded for $h\to 0$. Therefore, there
exists a sequence of positive numbers $h_n\; (n =1, 2, \ldots)$
such that $h_n \to 0$ as $n\to  \infty$, and the
sequence of functions $G_{h_n}(u)$ converges to a (bounded
nondecreasing) function $G(u)$ as $n\to  \infty$.
With
$$
\gamma_n = \int\limits_{-\infty}^{+\infty} \frac{d G_{h_n}(u)}{u}
$$
(where the integral has a sense due to (A.1)), we  have
$$
\log\, \varphi(t) = \lim_{n\to  \infty} \left\{i t \gamma_n
+  \int\limits_{-\infty}^{+\infty} \left(\e^{i t u} - 1 - \frac{i t
u}{1 + u^2}\right) \frac{u^2 + 1}{u^2}\,d G_{h_n}(u)\right\}\,.
$$
Since the  integrand of the above integral is bounded and continuous,
this integral tends as $n\to  \infty$ to  
$$
\int\limits_{-\infty}^{+\infty} \left(\e^{i t u} - 1 -
\frac{i t u}{1 + u^2}\right) \frac{u^2 + 1}{u^2}\,d G(u)\,.
$$
Hence the sequence $\gamma_n$ should converge to a certain
positive  constant $\gamma$. Therefore,
$$
\log\, \varphi(t) = i t \gamma + \int\limits_{-\infty}^{+\infty}
\left(e^{i t u} - 1 - \frac{i t u}{1 + u^2}\right) \frac{u^2 +
1}{u^2}\,d G(u)\,.
\eqno(A.5) $$
This  is the P.~L\'{e}vy formula up to certain unessential details
concerning the way of its presentation.

To prove the uniqueness of the last representation it is easier to
get the inversion formula. Let
$$
\Delta (t) = \int\limits_{t - 1}^{t + 1} \log\, \varphi(\alpha) d
\alpha - 2 \log\, [\varphi(t)]\,.
$$
Then  Eq. (A.5) gives  immediately
$$
\Delta (t) = - 2\int\limits_{-\infty}^{+\infty} \,e^{i t u}\left(1 -
\frac{\sin\, u}{u} \right)\, d G((u) =
\int\limits_{-\infty}^{+\infty} \e^{i t u} \,d K(u)\,,
$$
where
$$
K(u) = -2 \int\limits_{0}^{u} \left(1 - \frac{\sin\, v}{v} \right)
\,d G((v)\, .
$$
Then, the well-known P.~L\'{e}vy inversion formula\footnote{
{\it Calcul des probabilit\'{e}s}, Paris (1925),
p. 167. In the general case the integral should be considered in
the sense of the Cauchy principal value.}
yields
$$
K(u) = \frac{1}{2\pi} \int\limits_{-\infty}^{+\infty} \frac{1 -
\e^{-i t u}}{i t} \,\Delta(t) \,d t\,.
$$
It follows that $K(u)$ (and hence $G(u)$) is completely determined
by the function $\varphi(t)$. Then, we can easily conclude  that
$$
G_h(u) \to  G(u)\quad \hbox{as} \quad  h\to  0\,.
$$
If this were not true, then there it would exist a sequence of functions
$G_h(u)$ converging to another function (different from  $G(u)$).
This would  give another representation of the type (A.5) for the
function $\log\, \varphi(t)$.

Vice-versa, now we can show  that if the  logarithm of a function
$\varphi(t)$ is represented in the form (A.5) for a certain
nondecreasing bounded function $G(u)$, then $\varphi(t)$ is a
characteristic function of an infinitely divisible distribution.
Let $\varepsilon$ be an arbitrary positive number. Put
$$
\begin{array}{lll}
\Delta_{\varepsilon}&=& G(\varepsilon) - G(-\varepsilon), \\
&& \\
 G_{\varepsilon}(u)&=& \left\{
\begin{array}{ll}
G(u),& u\leq - \varepsilon,\\
G(-\varepsilon),& -\varepsilon \leq u \leq \varepsilon,\\
G(u) - \Delta_{\varepsilon},& u\geq \varepsilon.
\end{array}
\right.
\end{array}
$$
Since the function $G_{\varepsilon}(u)$ is bounded and
nondecreasing, we can write
$$
G_{\varepsilon}(u) = \lambda_\varepsilon \varphi_\varepsilon(u),
$$
where $\lambda_\varepsilon$ is a positive number and the function
$\varphi_\varepsilon(u)$ differs by an additive constant from a
certain distribution (this statement becomes trivial if the total
variation of $G_{\varepsilon}(u)$ is equal to zero for each
$\varepsilon > 0$). Let further
$$
f_\varepsilon(t) = \int\limits_{|u|>\varepsilon} \left(e^{i t u} -
1\right)d G(u) = \int\limits_{-\infty}^{+\infty} \left(e^{i t u} -
1\right)d G_\varepsilon(u) = $$
$$ = \lambda_\varepsilon
\int\limits_{-\infty}^{+\infty} \left(e^{i t u} - 1\right)d
\varphi_\varepsilon(u) = \lambda_\varepsilon
\left\{\varphi_\varepsilon(t) - 1\right\},
$$
where $\varphi_\varepsilon(t)$ is the characteristic function of
the distribution $\varphi_\varepsilon(x)$. Evidently, the expression
$$
\frac{\lambda_\varepsilon}{n} \varphi_\varepsilon(t) + \left\{1 -
\frac{\lambda_\varepsilon}{n}\right\}
$$
is a characteristic function for each $n\geq \lambda_\varepsilon$.
Hence the function
$$
n\log\, \left\{\frac{\lambda_\varepsilon}{n}
\varphi_\varepsilon(t) + \left(1 -
\frac{\lambda_\varepsilon}{n}\right)\right\} = n\log\, \left\{1 +
\frac{\lambda_\varepsilon}{n} \left[\varphi_\varepsilon(t) -
1\right]\right\}
$$
is the logarithm of a characteristic function. The same is true for
its limit as $n\to  \infty$ which is equal to
$$
\lambda_\varepsilon \left[\varphi_\varepsilon(t) - 1\right] =
f_\varepsilon(t)\,.
$$
Therefore, if $G(u)$ is an arbitrary bounded nondecreasing
function and $\varepsilon$ is an arbitrary positive number, then
the integral
$$
 \int\limits_{|u|>\varepsilon} \left(\e^{i t u} -
1\right)\,d G(u)
\eqno(A.6) $$
is the logarithm of a certain characteristic function. The same is
valid also for the integral
$$
 \int\limits_{|u|>\varepsilon} \left(\e^{i t u} -
1 - \frac{i t u}{1 + u^2}\right)\,d G(u)\,,
$$
which differs from (A.6) only by a  term $i t \gamma$, where
$\gamma$ is a real constant. We can also change in the last
integral $d G(u)$ to $(1 + u^2)/{u^2}\,d G(u)$, since the
function $(1 + u^2)/{u^2}$ is bounded for $|u| > \varepsilon$.
Finally we can pass to the limit as $\varepsilon \to  0$.
Therefore the function
$$
\lim_{\varepsilon \to  0} \int\limits_{|u|>\varepsilon}
\left(\e^{i t u} - 1 - \frac{i t u}{1 + u^2}\right)\frac{1 +
u^2}{u^2}\,d G(u)
$$
is the logarithm of a characteristic function. But the expression
(A.5) differs from this limit only by  a term $i t \gamma$
and a term of the type $- a t^2\; (a\geq 0)$ which is due to
a possible discontinuity of $G(u)$ at $u = 0$. Hence the
function $\varphi(t)$ is the product of a characteristic function
with an  expression of the type
$$
\e^{\, \ds i t \gamma - a t^2}\,,
$$
where $\gamma$ is a real constant and  $a \geq 0$. The last expression
is a characteristic function of a certain Gaussian Law. Hence the
function $\varphi(t)$ is a characteristic function as well.

The corresponding law is evidently infinitely divisible since
$\lambda \log\, \varphi(t)$ is for each $\lambda\geq 0$ the
expression of the same type as (A.5). Thus, by what is proved
above,
$\lambda \log\, \varphi(t)$ is the logarithm of a certain
characteristic function.

{\bf Supplement.} B.~V.~Gnedenko has pointed out that, to get the
statement   for the expression preceding to (A.5), one needs to
see that for $\alpha\to  \infty$ the 
limit
$$
\int\limits_{|u|\geq \alpha} d G_h(u) \to  0
$$
is {\it uniform with respect to $h$}.
To show this, it is  sufficient  to note
that, analogously to (A.4), one can prove the inequality
$$
-\frac{\alpha}{2} \int\limits_{0}^{2/\alpha} {\mathrm{Re}}
[I_h(t)]\, d t \geq \int\limits_{|u|\geq \alpha} \left(1 -
\frac{\sin\, (2u/ \alpha)}{(2u/\alpha)}\right) \,d G_h(u)
\geq {1\over 2} \int\limits_{|u|\geq \alpha}  d G_h(u).
$$
The left hand-side of this inequality tends as $h\to  0$
to
$$
-\frac{\alpha}{2} \int\limits_{0}^{2\over\alpha} {\mathrm{Re}}
[\log\, \varphi(t)]\, d t\,,
$$
which is sufficiently small for sufficiently large $\alpha$.



\end{document}